\documentclass[twosided,reqno]{amsart}

\usepackage{amsmath}
\usepackage{amssymb}
\usepackage{amsthm}

\theoremstyle{theorem}
\newtheorem{theorem}{\scshape Theorem }[section]

\newtheorem{proposition}[theorem]{\scshape Proposition}
\theoremstyle{definition}

\newcommand{\ma}{\mathbb}
\newcommand{\be}{\begin{equation}}
\newcommand{\ee}{\end{equation}}
\newcommand{\ben}{\begin{equation*}}
\newcommand{\een}{\end{equation*}}
\newcommand{\fa}{\frac}
\newcommand{\la}{\label}
\newcommand{\Z}{\mathbb{Z}}
\newcommand{\E}{E_{n}}
\newcommand{\eq}{E_{n,q}}
\newcommand{\ty}{\infty}
\newcommand{\I}{\int_{0}^{\infty}}
\newcommand{\U}{\sum_{n=0}^{\infty}}
\newcommand{\bt}{\begin{theorem}}
\newcommand{\et}{\end{theorem}}
\newcommand{\bi}{\binom}
\newcommand{\bp}{\begin{proposition}}
\newcommand{\ep}{\end{proposition}}

\numberwithin{equation}{section}

\begin{document}

\title{Symmetric identities of the $q$-Euler polynomials}

\author{Dae San Kim}
\address{Department of Mathematics, Sogang University, Seoul 121-742, Republic of Korea.}
\email{dskim@sogang.ac.kr}

\author{Taekyun Kim}
\address{Department of Mathematics, Kwangwoon University, Seoul 139-701, Republic of Korea}
\email{tkkim@kw.ac.kr}

\maketitle

\begin{abstract}
In this paper, we study some symmetric identities of $q$-Euler numbers and polynomials.  From these properties, we derive several identities of $q$-Euler numbers and polynomials.
\end{abstract}

\section{Introduction}

The Euler polynomials are defined by the generating function to be

\be\la{1}
\fa{2}{e^t+1}e^{xt}=e^{E(x)t}=\U \E(x)\fa{t^n}{n!},~(see~[2-5]).
\ee
with the usual convention about replacing $E^n(x)$ by $\E(x)$.\\

When $x=0$, $\E=\E(0)$ are called the Euler numbers. Let $q \in \ma{C}$ with $|q|<1$.  For any complex number $x$, the $q$-analogue of $x$ is defined by $[x]_q=\fa{1-q^x}{1-q}$. Note that $\lim_{q \to 1}[x]_q=x$. Recently, T. Kim introduced a $q$-extension of Euler polynomials as follows:

\be\la{2}
F_q(t,x)=[2]_q\U(-1)^nq^ne^{[n+x]_qt}=\U\eq(x)\fa{t^n}{n!},~(see[7,8]).
\ee

When $x=0$, $\eq=\eq(0)$ are called the $q$-Euler numbers.  From (\ref{2}), we note that

\be\la{3}
\begin{split}
\eq(x)&=(q^xE_q+[x]_q)^n\\
&=\sum_{l=0}^{n}\bi{n}{l}q^{xl}E_{l,q}[x]_q^{n-l},~(see~[7,8]),
\end{split}
\ee
with the usual convention about replacing $E_q^l$ by $E_{l,q}$.\\

In \cite{8}, Kim introduced $q$-Euler zeta function as follows:

\be\la{4}
\begin{split}
\zeta_{E,q}(s,x)&=\fa{1}{\Gamma(s)}\I t^{s-1}F_q(-t,x)dt\\
&=[2]_q \U \fa{(-1)^nq^n}{[n+x]_q^s},
\end{split}
\ee
where $x \neq 0,-1,-2,\ldots,$ and  $s \in \ma{C}$.\\

From (\ref{4}), we have
\be\la{5}
\zeta_{E,q}(-m,x)=E_{m,q}(x),
\ee
where $m \in \Z_{\geq 0}$.\\

Recently, Y. He gave some interesting symmetric identities of Carlitz's $q$-Bernoulli numbers and polynomials. In this paper, we study some new symmetries of the $q$-Euler numbers and polynomials, which is the answer to an open question for the symmetric identities of Carlitz's type $q$-Euler numbers and polynomials in \cite{5}.  By using our symmetries for the $q$-Euler polynomials we can obtain some identities between $q$-Euler numbers and polynomials.

\section{Symmetric identities of $q$-Euler polynomials}
In this section, we assume that $a,b \in \ma{N}$ with $a \equiv 1\ (\textrm{mod}\ 2)$ and $b \equiv1\ (\textrm{mod}\ 2)$.  First, we observe that
\be\la{6}
\begin{split}
&\fa{1}{[2]_{q^a}}\zeta_{E,q^a}(s,bx+\fa{bj}{a})=\U\fa{(-1)^nq^{na}}{[n+bx+\fa{bj}{a}]_{q^a}^s}\\
&=\U\fa{q^{an}(-1)^n[a]_q^s}{[bj+abx+an]_q^s}=[a]_q^s\U\sum_{i=0}^{b-1}\fa{(-1)^{i+bn}q^{a(i+bn)}}{[ab(x+n)+bj+ai]_q^s}.
\end{split}
\ee

Thus, by (\ref{6}), we get
\be\la{7}
\fa{[b]^s_q}{[2]_{q^a}}\sum_{j=0}^{a-1}(-1)^j q^{bj}\zeta_{E,q^a}(s,bx+\fa{bj}{a})
=[b]_q^s[a]_q^s\sum_{j=0}^{a-1}\sum_{i=0}^{b-1}\U \fa{q^{ai+bj+abn}(-1)^{i+n+j}}{[ab(x+n)+bj+ai]_q^s}.
\ee

By the same method as (\ref{7}), we get
\be\la{8}
\begin{split}
&\fa{[a]^s_q}{[2]_{q^b}}\sum_{j=0}^{b-1}(-1)^jq^{aj}\zeta_{E,q^b}(s,ax+\fa{aj}{b})\\
&=[a]_q^s[b]_q^s\sum_{j=0}^{b-1}\sum_{i=0}^{a-1}\U \fa{q^{bi+aj+abn}(-1)^{i+n+j}}{[ab(x+n)+aj+bi]_q^s}.
\end{split}
\ee

Therefore, by (\ref{7}) and (\ref{8}), we obtain the following theorem.
\bt\la{t1}
For $a,b\in \ma{N}$ with $a \equiv 1\ (\textrm{mod}\ 2)$, $b \equiv1\ (\textrm{mod}\ 2)$,
\ben
[2]_{q^b}[b]^s_q\sum_{j=0}^{a-1}(-1)^jq^{bj}\zeta_{E,q^a}(s,bx+\fa{bj}{a})=
[2]_{q^a}[a]^s_q\sum_{j=0}^{b-1}(-1)^jq^{aj}\zeta_{E,q^b}(s,ax+\fa{aj}{b}).
\een
\et

By (\ref{5}) and Theorem \ref{t1}, we obtain the following theorem.

\bt\la{t2}
For $n \in \Z_{\geq 0}$ and $a,b\in \ma{N}$ with $a \equiv 1\ (\textrm{mod}\ 2)$, $b \equiv1\ (\textrm{mod}\ 2)$, we have
\ben
[2]_{q^b}[a]^n_q\sum_{j=0}^{a-1}(-1)^jq^{bj}E_{n,q^a}(bx+\fa{bj}{a})=
[2]_{q^a}[b]^n_q\sum_{j=0}^{b-1}(-1)^jq^{aj}E_{n,q^b}(ax+\fa{aj}{b}).
\een
\et

From (\ref{3}), we note that

\be\la{9}
\begin{split}
E_{n,q}(x+y)&=(q^{x+y}E_q+[x+y]_q)^n\\
&=(q^{x+y}E_q+q^{x}[y]_q+[x]_q)^n\\
&=(q^x(q^yE_q+[y]_q)+[x]_q)^n\\
&=\sum_{i=0}^{n}\bi{n}{i}q^{xi}(q^yE_q+[y]_q)^i[x]_q^{n-i}\\
&=\sum_{i=0}^{n}\bi{n}{i}q^{xi}E_{i,q}(y)[x]_q^{n-i}.
\end{split}
\ee

Therefore, by (\ref{9}), we obtain the following proposition.

\bp\la{pro3}
For $n \geq 0$, we have
\ben
\begin{split}
E_{n,q}(x+y)&=\sum_{i=0}^{n}\bi{n}{i}q^{xi}E_{i,q}(y)[x]_q^{n-i}\\
&=\sum_{i=0}^{n}\bi{n}{i}q^{(n-i)x}E_{n-i,q}(y)[x]_q^{i}.
\end{split}
\een
\ep

Now, we observe that

\be\la{10}
\begin{split}
&\sum_{j=0}^{a-1}(-1)^jq^{bj}E_{n,q^a}(bx+\fa{bj}{a})\\
&=\sum_{j=0}^{a-1}(-1)^jq^{bj}\sum_{i=0}^{n}\bi{n}{i}q^{ia(\fa{bj}{a})}E_{i,q^a}(bx)\left[\fa{bj}{a}\right]_{q^a}^{n-i}\\
&=\sum_{j=0}^{a-1}(-1)^jq^{bj}\sum_{i=0}^{n}\bi{n}{i}q^{(n-i)bj}E_{n-i,q^a}(bx)\left[\fa{bj}{a}\right]_{q^a}^i\\
&=\sum_{i=0}^{n}\bi{n}{i}\left(\fa{[b]_q}{[a]_q}\right)^i E_{n-i,q^a}(bx)\sum_{j=0}^{a-1}(-1)^jq^{bj(n+1-i)}[j]^i_{q^b}\\
&=\sum_{i=0}^{n}\bi{n}{i}\left(\fa{[b]_q}{[a]_q}\right)^i E_{n-i,q^a}(bx)S_{n,i,q^b}^{*}(a),\\
\end{split}
\ee
where $S_{n,i,q}^{*}(a)=\sum_{j=0}^{a-1}(-1)^jq^{(n+1-i)j}[j]_q^i$.\\

From (\ref{10}), we can derive
\be\la{11}
[2]_{q^b}[a]_q^n\sum_{j=0}^{a-1}(-1)^jq^{bj}E_{n,q^a}(bx+\fa{bj}{a})=[2]_{q^b}\sum_{i=0}^{n}\bi{n}{i}[a]_q^{n-i}[b]_q^iE_{n-i,q^a}(bx)S_{n,i,q^b}^{*}(a).
\ee

By the same method as (\ref{11}), we get
\be\la{12}
[2]_{q^a}[b]_q^n\sum_{j=0}^{b-1}(-1)^jq^{aj}E_{n,q^b}(ax+\fa{aj}{b})=[2]_{q^a}\sum_{i=0}^{n}\bi{n}{i}[b]_q^{n-i}[a]_q^iE_{n-i,q^b}(ax)S_{n,i,q^a}^{*}(b).
\ee
Therefore, by Theorem 2.2, (\ref{11}) and (\ref{12}), we obtain the following theorem.

\bt\la{t4}
For $n \in \Z_{\geq 0}$ and $a,b\in \ma{N}$ with $a \equiv 1\ (\textrm{mod}\ 2)$, $b \equiv1\ (\textrm{mod}\ 2)$, we have
\ben
[2]_{q^b}\sum_{i=0}^{n}\bi{n}{i}[a]_q^{n-i}[b]_q^iE_{n-i,q^a}(bx)S_{n,i,q^b}^{*}(a)=[2]_{q^a}\sum_{i=0}^{n}\bi{n}{i}[b]_q^{n-i}[a]_q^iE_{n-i,q^b}(ax)S_{n,i,q^a}^{*}(b),
\een
where $S_{n,i,q}^{*}(a)=\sum_{j=0}^{a-1}(-1)^jq^{(n+1-i)j}[j]_q^i$.
\et

It is easy to show that
\be\la{13}
[x]_{q}u+q^x[y+m]_q(u+v)=[x+y+m]_q(u+v)-[x]_qv.
\ee

Thus, by (\ref{13}), we get
\be\la{14}
e^{[x]_q u}\sum_{m=0}^{\ty}q^m (-1)^m e^{[y+m]_q q^x(u+v)}=e^{-[x]_qv}\sum_{m=0}^{\ty}q^m (-1)^mq^{[x+y+m]_q(u+v)}.
\ee

The left hand side of (\ref{14}) multiplied by $[2]_q$ is given by

\be\la{15}
\begin{split}
&[2]_qe^{[x]_qu}\sum_{m=0}^{\ty}q^m (-1)^m e^{[y+m]_q q^x (u+v)}\\
&=e^{[x]_qu}\U q^{nx}E_{n,q}(y) \fa{(u+v)^n}{n!}\\
&=\left(\sum_{l=0}^{\ty}[x]_q^l\fa{u^l}{l!} \right) \left( \sum_{k=0}^{\ty} \U q^{(k+n)x}E_{k+n,q}(y)\fa{u^k}{k!}\fa{v^n}{n!}\right)\\
&=\sum_{m=0}^{\ty}\U \left( \sum_{k=0}^{m}\bi{m}{k}q^{(k+n)x}E_{k+n,q}(y)[x]_q^{m-k}\right) \fa{u^m}{m!}\fa{v^n}{n!}.
\end{split}
\ee

The right hand side of (\ref{14}) multiplied by $[2]_q$ is given by
\be\la{16}
\begin{split}
&[2]_qe^{-[x]_qv}\sum_{m=0}^{\ty}(-1)^m q^m e^{[x+y+m]_q(u+v)}\\
&=e^{-[x]_qv} \U \eq(x+y)\fa{(u+v)^n}{n!}\\
&=\left(\sum_{l=0}^{\ty}\fa{(-[x]_q)^l}{l!}v^l\right)\left(\sum_{m=0}^{\ty}\sum_{k=0}^{\ty}E_{m+k,q}(x+y)\fa{u^m}{m!}\fa{v^k}{k!}\right)\\
&=\U\sum_{m=0}^{\ty}\left( \sum_{k=0}^{n}\bi{n}{k}E_{m+k,q}(x+y)(-[x]_q)^{n-k}\right)\fa{u^m}{m!}\fa{v^n}{n!}\\
&=\U\sum_{m=0}^{\ty}\left( \sum_{k=0}^{n}\bi{n}{k}E_{m+k,q}(x+y)q^{(n-k)x}[-x]_q^{n-k}\right)\fa{u^m}{m!}\fa{v^n}{n!}.
\end{split}
\ee
Therefore, by (\ref{15}) and (\ref{16}), we get
\be\la{17}
 \sum_{k=0}^{m}\bi{m}{k}q^{(n+k)x}E_{n+k,q}(y)[x]_q^{m-k}= \sum_{k=0}^{n}\bi{n}{k}q^{(n-k)x}E_{m+k,q}(x+y)[-x]_q^{n-k}
\ee



\begin{thebibliography}{99}
\bibitem{1}
Y. He,
\emph {Symmetric identities for Calitz's $q$-Bernoulli numbers and polynomials}, Advances in Difference Equations 2013(2013), 246, doi:10.1186/1687-1847-2013-246.

\bibitem{2}
D. S. Kim,
\emph{Symmetry identities for generalized twisted Euler polynomials twisted by unramified roots of unity}, Proc. Jangjeon Math. Soc. 15(2012), no. 3, 303-316.

\bibitem {3} D. S. Kim, N. Lee, J. Na and K. H. Park, {\it Identities of symmetry for higher-order Euler polynomials in three variables (I)}, Adv. Stud. Contemp. Math.  ${\mathbf{22}}$  (2012),  no. 1, 51-74.

\bibitem {4} D. S. Kim, N. Lee, J. Na and K. H. Park, {\it Identities of symmetry for higher-order Euler polynomials in three variables (II)}, J. Math. Anal. Appl.  ${\mathbf{379}}$  (2011),  no. 1, 388-400.

\bibitem{5}
T. Kim,
\emph{Symmetry $p$-adic invariant integral on $\Z_p$ for Bernoulli and Euler polynomials}, J. Difference Equ. Appl. 14(2008), no.12, 1267-1277.

\bibitem {6} T. Kim, {\it An identity of the symmetry for the Frobenius-Euler polynomials associated with the fermionic $p$-adic invariant q-integrals on $\Bbb Z_p$}, Rocky Mountain J. Math. ${\mathbf{41}}$ (2011),  no. 1, 239-247.

\bibitem{7}
T. Kim,
\emph{Analytic continuation of $q$-Euler numbers and polynomials},
Applied mathematics Letters 21(2008), 1320-1323.

\bibitem{8}
T. Kim,
\emph{$q$-Euler numbers and polynomials associated with $p$-adic $q$-integrals}, Journal of Nonlinear Mathematical Physics, 14(2007), no.1, 15-27.
\end{thebibliography}
\end{document}